\documentclass[10pt,twocolumn,a4paper]{article}

\usepackage[left=20mm, right=15mm, top=15mm, bottom=15mm]{geometry}

 \usepackage{amsthm}
\usepackage{amsfonts}
\usepackage{bbm}

\usepackage{subfig}
\usepackage{flushend}
\usepackage{indentfirst}
\usepackage{graphics}
\usepackage{amsmath}
\usepackage{graphicx}
\usepackage{epstopdf}
\usepackage{float}

\usepackage[font=footnotesize,labelfont=bf]{caption}
\usepackage[labelsep=period]{caption}

\graphicspath{{./figure/}}

\newtheorem{ozn}{Definition}
\newtheorem{thm}{Theorem}

\begin{document}

\title{\huge \textbf{Minimax Estimation Problem for Periodically Correlated Stochastic Processes}}

\date{}

\twocolumn[

\begin{@twocolumnfalse}
\maketitle

\noindent Journal of Mathematics and System Science  Vol. 3, No.1,  26-30, 2013\\
ISSN 2159-5291, USA

\vspace{10pt}

\author{\textbf{Iryna Dubovets'ka}, \textbf{Mykhailo Moklyachuk}$^{*}$, \\\\
 {Department of Probability Theory, Statistics and Actuarial
Mathematics, \\
Taras Shevchenko National University of Kyiv, Kyiv 01601, Ukraine}\\
$^{*}$Corresponding Author: Moklyachuk@gmail.com}\\\\\\

\end{@twocolumnfalse}

]

\noindent \textbf{\Large{Abstract}} \hspace{2pt}
The problem of optimal linear estimation of the functional $A\zeta =\int_{0}^{\infty }{a}(t)\zeta (t)dt$ depending on the unknown values of periodically correlated stochastic process
from observations of this process $\zeta (t)$ at points $t<0$  is considered.
 Formulas that determine the greatest value of the mean square error and the minimax estimation for the functional are proposed for a given class of admissible processes.
\\

\noindent \textbf{\Large{Keywords}} \hspace{2pt}  Periodically correlated process, minimax estimation, mean square error, least favorable process.\\

\section{\Large{Introduction}}

Methods of solution of problems of estimation of unknown values of stationary stochastic processes
 (extrapolation, interpolation and filtering problems) are developed in the works of Kolmogorov [1], Wiener [2], Yaglom [3, 4].
   These methods are based on the assumption that spectral densities of processes are exactly known.
   In the case where complete information on the spectral densities is impossible, while a set of admissible densities is given, the minimax approach to estimation problem is reasonable.
   That is we find the estimate that minimizes the mean square error for all spectral densities from a given class simultaneously.
   Ulf Grenander [5] was the first who applied the minimax estimation method to the extrapolation problem for stationary processes.
   Moklyachuk [6], Moklyachuk and Masyutka [7] studied the extrapolation, interpolation and filtering problems for stationary processes and sequences.
   In the paper by Gladyshev [8] the investigation of periodically correlated processes was started.
   The analysis of properties of correlation function and representation of periodically correlated processes is presented.
   Relations between periodically correlated processes and stationary processes were investigated by Makagon [9, 10].
   Minimax estimation problems for
   linear functionals from periodically correlated sequences were studied in works by Dubovets'ka,
   Masyutka and Moklyachuk [11-13].

\section{\Large{Periodically correlated processes and generated vector stationary sequences}}

\begin{ozn} $[8]$
 A mean square continuous stochastic process $\zeta :\mathbb{R}\to H={{L}_{2}}(\Omega ,\mathcal{F},\mathbb{P})$, $E\zeta (t)=0,$ is called periodically correlated (PC) with period $T$, if its correlation function $K(t+u,u)=E\zeta (t+u)\overline{\zeta (u)}$ for all $t,u\in \mathbb{R}$ and some fixed $T$is such that
\[K(t+u,u)=K(t+u+T,u+T).\]
\end{ozn}

  Consider the problem of optimal linear estimation of the functional
\[A\zeta =\int_{0}^{\infty }{a}(t)\zeta (t)dt\]
depending on the unknown values of PC stochastic process $\zeta (t)$ from the class $\mathbf{Y}$ of mean square continuous PC processes $\zeta (t)$ such that
$E\zeta (t)=0$, $E|\zeta (t){{|}^{2}}\le P/T$.
The estimation is based on observations of the process $\zeta (t)$ for $t<0$.
The function $a(t),$ $t\in \mathbb{R}$ satisfies the condition
$\int_{0}^{\infty }{|}a(t)|dt<\infty $.

   The functional $A\zeta $ can be represented in the form
\[A\zeta =\int_{0}^{\infty }{a}(t)\zeta (t)dt=\sum\nolimits_{j=0}^{\infty }{\int_{0}^{T}{{{a}_{j}}(u){{\zeta }_{j}}(u)du}},\]
\[a(u+jT)={{a}_{j}}(u),\zeta (u+jT)={{\zeta }_{j}}(u),u\in [0,T).\]
The sequence $\{{{\zeta }_{j}}={{\zeta }_{j}}(u),\,u\in [0,T),\,\,j\in \mathbb{Z}\}$ is ${{L}_{2}}([0,T);H)$-valued stationary stochastic sequence  with the correlation function
\[B(l,j)={{\langle {{\zeta }_{l}},{{\zeta }_{j}}\rangle }_{H}}=\int_{0}^{T}{K}(u+(l-j)T,u)du=B(l-j).\]

Define in the space ${{L}_{2}}([0,T);\mathbb{R})$ the orthonormal basis
\[\left\{{\tilde{e}}_{k}={T}^{-\frac{1}{2}}{{e}^{2\pi i\left\{{{(-1)}^{k}}\left[ {\scriptstyle{}^{k}/{}_{2}} \right]\right\}u/T}},k=1,2,3,\ldots \right\},\]
 \[\langle {{\tilde{e}}_{j}},{{\tilde{e}}_{k}}\rangle ={{\delta }_{kj}}.\]
With the help of this basis the stationary sequence $\{{{\zeta }_{j}},j\in \mathbb{Z}\}$ can be represented in the form
\[{{\zeta }_{j}}=\sum\nolimits_{k=1}^{\infty }{{{\zeta }_{kj}}}{{\tilde{e}}_{k}},\quad \]\[{{\zeta }_{kj}}=\langle {{\zeta }_{j}},{{\tilde{e}}_{k}}\rangle ={{T}^{-{\scriptstyle{}^{1}/{}_{2}}}}\int_{0}^{T}{{{\zeta }_{j}}}(v){{e}^{-2\pi i\{{{(-1)}^{k}}\left[ {\scriptstyle{}^{k}/{}_{2}} \right]\}v/T}}dv,\]
and the functional $A\zeta$ can be represented in the form
\[A\zeta =\sum\nolimits_{j=0}^{\infty }{\sum\nolimits_{k=1}^{\infty }{{{a}_{kj}}{{\zeta }_{kj}}}}=\sum\nolimits_{j=0}^{\infty }{\vec{a}_{j}^{\top }{{{\vec{\zeta }}}_{j}}},\]
\[{{\vec{\zeta }}_{j}}=({{\zeta }_{kj}},k=1,2,\ldots )^{\top },\] \[{{\vec{a}}_{j}}=({{a}_{kj}},k=1,2,\ldots )^{\top }=\]
\[={{({{a}_{1j}},{{a}_{3j}},{{a}_{2j}},\ldots ,{{a}_{2k+1,j}},{{a}_{2k,j}},\ldots )}^{\top }},\quad {{a}_{kj}}=\langle {{a}_{j}},{{\tilde{e}}_{k}}\rangle .\]
Components ${{\zeta }_{kj}}$ of the stationary sequence $\{{{\zeta }_{j}},j\in \mathbb{Z}\}$ are such that
\[E{{\zeta }_{kj}}=0,\]
\[\|{{\zeta }_{j}}\|_{H}^{2}=\sum\nolimits_{k=1}^{\infty }{E|{{\zeta }_{kj}}{{|}^{2}}}\le P,\]
\[E{{\zeta }_{kl}}\overline{{{\zeta }_{nj}}}=\langle R(l-j){{e}_{k}},{{e}_{n}}\rangle, \]
where $\{e_k,k=1,2,\dots\}$ is a basis in the space in ${{\ell }_{2}}$.
\\
The correlation function $R(n)$ of stationary sequence $\{{{\zeta }_{j}},j\in \mathbb{Z}\}$ is a correlation operator function in ${{\ell }_{2}}$.

\begin{ozn}
Denote by $H_{\zeta}(n)$ closed linear subspace of the Hilbert space $H$ generated by random values $\{{\zeta }_{kj},k\geq 1,j\leq n\}$
The sequence of stationary sequence $\{{{\zeta }_{j}},j\in \mathbb{Z}\}$ is called regular if $\cap_{n}H_{\zeta}(n)=\emptyset$.
In the case  $\cap_{n}H_{\zeta}(n)=H$ the sequence $\{{{\zeta }_{j}},j\in \mathbb{Z}\}$  is called singular.
\end{ozn}
 A stationary sequence $\{{{\zeta }_{j}},\,j\in \mathbb{Z}\}$ admits a unique representation in the form
${{\zeta }_{j}}=\zeta _{j}^{r}+\zeta _{j}^{s}$
where $\{\zeta _{j}^{r},\,j\in \mathbb{Z}\}$ is a regular sequence and $\{\zeta _{j}^{s},\,j\in \mathbb{Z}\}$ is a singular sequence.
Moreover, the sequences $\{\zeta _{j}^{r}\}$ and $\{\zeta _{n}^{s}\}$ are orthogonal for all $j,n\in \mathbb{Z}$.

Since the unknown values of components of singular stationary sequence have error-free estimate, we will consider the estimation problem only for regular
stationary sequences.

The regular stationary sequence $\{{{\zeta }_{j}},j\in \mathbb{Z}\}$ has the spectral density operator function $f(\lambda ),$ $\lambda \in [-\pi ,\pi )$ in ${{\ell }_{2}}$ and satisfy the equality
\[\langle R(l-j){{e}_{k}},{{e}_{n}}\rangle =\frac{1}{2\pi }\int_{-\pi }^{\pi }{{{e}^{i(l-j)\lambda }}}\langle f(\lambda ){{e}_{k}},{{e}_{n}}\rangle d\lambda ,\, k,n\ge 1.\]
The spectral density $f(\lambda )$ a.e. on $[-\pi ,\pi )$ is a kernel operator with integrable kernel norm
\begin{equation}\label{eq1}
\sum\nolimits_{k=1}^{\infty }{\frac{1}{2\pi }}\int_{-\pi }^{\pi }{\langle f(}\lambda ){{e}_{k}},{{e}_{k}}\rangle d\lambda =\|{{\zeta }_{j}}\|_{H}^{2}\le P.
\end{equation}
The regular stationary sequence $\{{{\zeta }_{j}},j\in \mathbb{Z}\}$ admits the canonical moving average representation of components [14]
 \begin{equation}\label{eq2}
{{\zeta }_{kj}}=\sum\nolimits_{s=-\infty }^{j}{\sum\nolimits_{m=1}^{M}{{{g}_{km}}(j-s){{\varepsilon }_{m}}(s)}},
\end{equation}
where ${{\varepsilon }_{m}}(s),$ $m=1,\ldots ,M,$ $s\in \mathbb{Z}$ are mutually orthogonal sequences in $H$ with orthonormal values; $M$ is the multiplicity of
$\{{{\zeta }_{j}}\}$; and ${{g}_{km}}(s),$ $k=1,2,\ldots ,$ $m=1,\ldots ,M,s=0,1,\ldots,$ are real-valued sequences such that
$$\sum\nolimits_{s=0}^{\infty }{\sum\nolimits_{k=1}^{\infty }{\sum\nolimits_{m=1}^{M}{{{\left| {{g}_{km}}(s) \right|}^{2}}}}}\le P.$$
As a consequence of the representation \eqref{eq2} the optimal linear estimation of components of stationary sequence $\{{{\zeta }_{j}},j\in \mathbb{Z}\}$ can be written in the form
 \begin{equation}\label{eq3}
{{\hat{\zeta }}_{kj}}=\sum\nolimits_{s=-\infty }^{-1}{\sum\nolimits_{m=1}^{M}{{{g}_{km}}(j-s){{\varepsilon }_{m}}(s)}}.
\end{equation}
   The spectral density $f(\lambda )$ of regular stationary sequence $\{{{\zeta }_{j}},j\in \mathbb{Z}\}$ admits canonical factorization
    \begin{equation}\label{eq4}
{{f}_{kn}}(\lambda )=\sum\nolimits_{m=1}^{M}{{{{\hat{g}}}_{km}}({{e}^{i\lambda }})\overline{{{{\hat{g}}}_{nm}}({{e}^{i\lambda }})}},\quad k,n\ge 1,
\end{equation}
\[{{\hat{g}}_{km}}({{e}^{i\lambda }})=\sum\nolimits_{s=0}^{\infty }{{{g}_{km}}(s){{e}^{-is\lambda }}}.\]
This means that $$f(\lambda )=G(\lambda ){{G}^{*}}(\lambda ),\quad G(\lambda )=\{{{\hat{g}}_{km}}({{e}^{i\lambda }})\}_{k=\overline{1,\infty }}^{m=\overline{1,M}}.$$

\section{\Large{Minimax estimation of linear functional}}

   Let us assume that coefficients $\{{{\vec{a}}_{j}},j\ge 0\}$ satisfy conditions
  \begin{equation}\label{eq5}
	\sum\nolimits_{j=0}^{\infty }{\|{{{\vec{a}}}_{j}}\|}<\infty ,\quad \sum\nolimits_{j=0}^{\infty }{(j+1)\|{{{\vec{a}}}_{j}}{{\|}^{2}}}<\infty,
\end{equation}
\[\|{{\vec{a}}_{j}}{{\|}^{2}}=\sum\nolimits_{k=1}^{\infty }{|{{a}_{kj}}{{|}^{2}}}.\]
Denote by $\Lambda $ the set of all linear estimates of the functional $A\zeta $ based on observations of the process $\zeta (t)$ for $t<0$.
Let ${{\mathbf{Y}}_{R}}$ denotes the class of all regular stationary sequences that satisfy the condition $\|{{\zeta }_{j}}\|_{H}^{2}\le P$.

   We calculate the largest values of mean square error $\Delta (\zeta ,\hat{A})=E|A\zeta -\hat{A}\zeta {{|}^{2}}$ of estimation $\hat{A}\zeta $ of the functional $A\zeta $ and
   the largest values of mean square error $\Delta (\zeta ,{{\hat{A}}_{N}})=E|{{A}_{N}}\zeta -{{\hat{A}}_{N}}\zeta {{|}^{2}}$ of estimation ${{\hat{A}}_{N}}\zeta $ of the functional ${{A}_{N}}\zeta =\sum\nolimits_{j=0}^{N}{\vec{a}_{j}^{\top }{{{\vec{\zeta }}}_{j}}}$.

\begin{thm} \label{thm3.1}
 Let the coefficients $\{{{\vec{a}}_{j}},j=0,...,N\}$ which determine the functional ${{A}_{N}}\zeta$
 satisfy the condition
 \begin{equation}\label{eq6}
	\|{{\vec{a}}_{j}}\|<\infty ,\quad j=0,1,\ldots ,N.
\end{equation}
The function $\Delta (\zeta ,{{\hat{A}}_{N}})$ has a saddle point on the set $\mathbf{Y}\times \Lambda $ and
\[\underset{{{{\hat{A}}}_{N}}\in \Lambda }{\mathop{\min }}\,\underset{\zeta \in \mathbf{Y}}{\mathop{\max }}\,\Delta (\zeta ,{{\hat{A}}_{N}})=\underset{\zeta \in \mathbf{Y}}{\mathop{\max }}\,\underset{{{{\hat{A}}}_{N}}\in \Lambda }{\mathop{\min }}\,\Delta (\zeta ,{{\hat{A}}_{N}})=P\cdot \nu _{N}^{2},\]
where $\nu _{N}^{2}$ is the greatest eigenvalue of the self-adjoint compact operator ${{Q}_{N}}=\{{{Q}_{N}}(p,q)\}_{p,q=0}^{N}$ in the space ${{\ell }_{2}}$ determined by block-matrices ${{Q}_{N}}(p,q)=\{Q_{kn}^{N}(p,q)\}_{k,n=1}^{\infty }$ with elements
 \begin{equation}\label{eq7}
Q_{kn}^{N}(p,q)=\sum\nolimits_{s=0}^{min(N-p,N-q)}{{{a}_{k,s+p}}\cdot \overline{{{a}_{n,s+q}}}},
\end{equation}
\[k,n=1,2,3,\ldots ,\quad p,q=0,1,\ldots ,N.\]
The least favorable stochastic sequence in the class $\mathbf{Y}$ for the optimal estimate of the functional $A_N\zeta $ is one-sided moving average sequence of order $N$ of the form
\[{{\vec{\zeta }}_{j}}=\sum\nolimits_{s=j-N}^{j}{g(j-s)\vec{\varepsilon }(s)},\]
where $g=(g(p))_{p=0}^{N}$ is the eigenvector, that corresponds to $\nu_N^2$, constructed from matrices $g(s)=\{{{g}_{km}}(s)\}_{k=\overline{1,\infty }}^{m=\overline{1,M}}$  and determined by the condition $\|g{{\|}^{2}}=\sum\nolimits_{p=0}^{N }{\|g(p){{\|}^{2}}}=P$, and $\vec{\varepsilon }(s)=\{{{\varepsilon }_{m}}(s)\}_{m=1}^{M}$ is a vector stationary stochastic sequence with orthogonal values.
 \end{thm}

   Proof: Upper bound. Denote by ${{\Lambda }_{1}}$ the class of all linear estimates of the functional ${{A}_{N}}\zeta $ of the form ${{\hat{A}}_{N}}\zeta =\sum\nolimits_{j=-\infty }^{-1}{\vec{c}_{j}^{\top }{{{\vec{\zeta }}}_{j}}}.$ Taking into account the spectral decomposition of stationary sequence and of its correlation function [15], the following relation holds true
\[\Delta (\zeta ,{{\hat{A}}_{N}})=\]
\[=\frac{1}{2\pi }\int_{-\pi }^{\pi }{{{\left( {{A}_{N}}({{e}^{i\lambda }})-C({{e}^{i\lambda }}) \right)}^{*}}}f(\lambda )\left( {{A}_{N}}({{e}^{i\lambda }})-C({{e}^{i\lambda }}) \right)d\lambda ,\]
	\[{{A}_{N}}({{e}^{i\lambda }})=\sum\nolimits_{j=0}^{N }{{{{\vec{a}}}_{j}}{{e}^{ij\lambda }}},\quad C({{e}^{i\lambda }})=\sum\nolimits_{j=-\infty }^{-1}{{{{\vec{c}}}_{j}}{{e}^{ij\lambda }}}.
\]
For the mean square error, we can derive the following estimates
\[\underset{\zeta \in \mathbf{Y}}{\mathop{\max }}\,\Delta (\zeta ,{{\hat{A}}_{N}})\le
\]
\[\underset{\zeta \in \mathbf{Y}}{\mathop{\max }}\,\frac{1}{2\pi }\int_{-\pi }^{\pi }{\|{{A}_{N}}(}{{e}^{i\lambda }})-C({{e}^{i\lambda }}){{\|}^{2}}\|f(\lambda )\|d\lambda \le \]\[\le \underset{\lambda \in [-\pi ,\pi )}{\mathop{\max }}\,\|{{A}_{N}}({{e}^{i\lambda }})-C({{e}^{i\lambda }}){{\|}^{2}}\cdot \frac{1}{2\pi }\int_{-\pi }^{\pi }{\|f(}\lambda )\|d\lambda \le \]
\[\le P\cdot \underset{\lambda \in [-\pi ,\pi )}{\mathop{\max }}\,\|{{A}_{N}}({{e}^{i\lambda }})-C({{e}^{i\lambda }}){{\|}^{2}},\]
what follows from the inequality  $$\|f(\lambda )\|={{\left( \sum\nolimits_{k,n=1}^{\infty }{|{{f}_{kn}}(\lambda ){{|}^{2}}} \right)}^{1/2}}\le \text{Tr}\,f(\lambda )$$ and integrability of  kernel norm of spectral density \eqref{eq1}.

   To calculate $\underset{\lambda \in [-\pi ,\pi )}{\mathop{\max }}\,\|{{A}_{N}}({{e}^{i\lambda }})-C({{e}^{i\lambda }}){{\|}^{2}}$
   we consider the class of all power series $\vec{f}(z)=\sum\nolimits_{j=0}^{\infty }{{{{\vec{\alpha }}}_{j}}{{z}^{j}}}$, which are regular for $|z|<1$  and have first $(N+1)$ summands $\sum\nolimits_{j=0}^{N}{{{{\vec{d}}}_{j}}{{z}^{j}}}$ fixed. Denote by $\omega _{N}^{2}$ the greatest eigenvalue of the matrix  ${{D}_{N}}=\{{{D}_{N}}(p,q)\}_{p,q=0}^{N}$ constructed by block-matrices
\[{{D}_{N}}(p,q)=\sum\nolimits_{s=0}^{min(N-p,N-q)}{{{{\vec{a}}}_{N-p+s}}}\cdot \,\vec{a}_{N-q+s}^{*},\] \[ p,q=0,1,\ldots ,N.\]
Then $\underset{\{{{{\vec{\alpha }}}_{j}}:j\ge N+1\}}{\mathop{\min }}\,\underset{|z|=1}{\mathop{\max }}\,\|\vec{f}(z){{\|}^{2}}=\omega _{N}^{2}$
as it follows from [16].
Since ${{\Lambda }_{1}}\subset \Lambda $, we have
\begin{equation}\label{eq8}
\underset{{{{\hat{A}}}_{N}}\in \Lambda }{\mathop{\min }}\,\underset{\zeta \in \mathbf{Y}}{\mathop{\max }}\,\Delta (\zeta ,{{\hat{A}}_{N}})\le \underset{{{{\hat{A}}}_{N}}\in {{\Lambda }_{1}}}{\mathop{\min }}\,\underset{\zeta \in \mathbf{Y}}{\mathop{\max }}\,\Delta (\zeta ,{{\hat{A}}_{N}})\le P\cdot \omega _{N}^{2}.
\end{equation}

   Lower bound. Using the canonical decompositions of components of regular stationary sequence \eqref{eq2} and components of optimal linear estimation \eqref{eq3}, we obtain
\[\underset{\hat{A}\in \Lambda }{\mathop{\min }}\,\Delta (\zeta ,\hat{A}_N)=E{{\left| \sum\nolimits_{k=1}^{\infty }{\sum\nolimits_{j=0}^{N}{{{a}_{kj}}}({{\zeta }_{kj}}-{{{\hat{\zeta }}}_{kj}})} \right|}^{2}}=\]
\begin{equation}\label{eq9}
=\sum\nolimits_{k,n=1}^{\infty }{\sum\nolimits_{m=1}^{M}{\sum\nolimits_{p,q=0}^{N}{{{g}_{km}}(p)\overline{{{g}_{nm}}(q)}Q_{kn}^{N}(p,q)}}},
\end{equation}
where operator ${{Q}_{N}}=\{{{Q}_{N}}(p,q)\}_{p,q=0}^{N}$ in the space ${{\ell }_{2}}$ determined by block-matrices ${{Q}_{N}}(p,q)=\{Q_{kn}^{N}(p,q)\}_{k,n=1}^{\infty }$
with elements \eqref{eq7}.
It can be represented in the form ${{Q}_{N}}={{A}_{N}}A_{N}^{*}$, where the matrix operator ${{A}_{N}}=\{{{A}_{N}}(p,q)\}_{p,q=0}^{\infty }$
is determined by the vector columns ${{A}_{N}}(p,q)={{\vec{a}}_{p+q}},p+q\le N,$ and ${{A}_{N}}(p,q)=\vec{0},p+q>N.$
The squared absolute norm of operator ${{A}_{N}}$ equals
${{\mathcal{N}}^{2}}({{A}_{N}})=$ $\sum\nolimits_{p=0}^{\infty }{(p+1)}\|{{\vec{a}}_{p}}{{\|}^{2}}<\infty $, because of \eqref{eq6}.
Then operator ${{A}_{N}}$ is continuous [17]. The operator ${{Q}_{N}}$ is continuous and has real nonnegative eigenvalues.

   Denote by ${{g}_{N}}={{(g(0),g(1),\ldots ,g(N))}^{\top }}$ the vector columns with matrix elements$$g(p)=\{{{g}_{km}}(p)\}_{k=\overline{1,\infty }}^{m=\overline{1,M}},p=0,1,\ldots ,N.$$
    Operator ${{Q}_{N}}$ acts on the vector ${{g}_{N}}$ as follows $${{Q}_{N}}{{g}_{N}}=\left\{\sum\nolimits_{q=0}^{N}{{{Q}_{N}}(p,q)g(q)}\right\}_{p=0}^{N}.$$
    The relation \eqref{eq9} has the form
\[\underset{{{{\hat{A}}}_{N}}\in \Lambda }{\mathop{\min }}\,\Delta (\zeta ,{{\hat{A}}_{N}})=\sum\nolimits_{p=0}^{N}{\sum\nolimits_{q=0}^{N}{\left({{Q}_{N}}(p,q)\overline{g(q)},\overline{g(p)}\right)}}=\]
\[=\left( {{Q}_{N}}\overline{{{g}_{N}}},\overline{{{g}_{N}}} \right),\]
where $(\cdot,\cdot)$ is the scalar product in $\ell_2$.

Let us denote  $\tilde{g}=\overline{{{g}_{N}}}{{P}^{-{\scriptstyle{}^{1}/{}_{2}}}}$.
Then the extremal problem \[|({{Q}_{N}}\tilde{g},\tilde{g})|\to \max,\,\,\|\tilde{g}\|=1\]
has solutions [17]. This solution is the eigenvalue of the operator ${{Q}_{N}}$.
The corresponding eigenvalue is $$\nu _{N}^{2}=\underset{\|\tilde{g}\|\le 1}{\mathop{\max }}\,|({{Q}_{N}}\tilde{g},\tilde{g})|=\|{{Q}_{N}}\|.$$
Thus we have
\[\underset{\zeta \in {{\mathbf{Y}}_{\mathbf{R}}}}{\mathop{\max }}\,\underset{{{{\hat{A}}}_{N}}\in \Lambda }{\mathop{\min }}\,\Delta (\zeta ,{{\hat{A}}_{N}})=P\cdot \underset{\|\tilde{g}\|\le 1}{\mathop{\max }}\,|({{Q}_{N}}\tilde{g},\tilde{g})|=\]
\begin{equation}\label{eq10}
=P\cdot \nu _{N}^{2}\le \underset{\zeta \in \mathbf{Y}}{\mathop{\max }}\,\underset{{{{\hat{A}}}_{N}}\in \Lambda }{\mathop{\min }}\,\Delta (\zeta ,{{\hat{A}}_{N}}).
\end{equation}

   Let us note that ${{D}_{N}}(N-p,N-q)={{Q}_{N}}(p,q)$. For this reason $\omega _{N}^{2}=\nu _{N}^{2}$. Relations \eqref{eq8} and \eqref{eq10} give us the inequality
   \begin{equation}\label{eq11}
	\underset{{{{\hat{A}}}_{N}}\in \Lambda }{\mathop{\min }}\,\underset{\zeta \in \mathbf{Y}}{\mathop{\max }}\,\Delta (\zeta ,{{\hat{A}}_{N}})\le P\cdot \nu _{N}^{2}\le \underset{\zeta \in \mathbf{Y}}{\mathop{\max }}\,\underset{{{{\hat{A}}}_{N}}\in \Lambda }{\mathop{\min }}\,\Delta (\zeta ,{{\hat{A}}_{N}}).
\end{equation}

Since the opposite inequality always holds true, we have equality in \eqref{eq11}. Theorem is proved.

   \begin{thm} \label{thm3.2}
Let the coefficients $\{{{\vec{a}}_{j}},j=0,1,...\}$ satisfy conditions \eqref{eq5}. Then the function $\Delta (\zeta ,\hat{A})$ has a saddle point on the set $\mathbf{Y}\times \Lambda $ and the following equality holds true
\begin{equation}\label{eq12}
\underset{\hat{A}\in \Lambda }{\mathop{\min }}\,\underset{\zeta \in \mathbf{Y}}{\mathop{\max }}\,\Delta (\zeta ,\hat{A})=\underset{\zeta \in \mathbf{Y}}{\mathop{\max }}\,\underset{\hat{A}\in \Lambda }{\mathop{\min }}\,\Delta (\zeta ,\hat{A})=P\cdot {{\nu }^{2}},
\end{equation}
where ${{\nu }^{2}}$ is the greatest eigenvalue of the self-adjoint compact operator $Q=\{Q(p,q)\}_{p,q=0}^{\infty }$ in the space ${{\ell }_{2}}$ determined by block-matrices $Q(p,q)=\{Q_{kn}^{{}}(p,q)\}_{k,n=1}^{\infty }$ with elements
\[Q_{kn}^{{}}(p,q)=\sum\nolimits_{s=0}^{\infty }{{{a}_{k,s+p}}\cdot \overline{{{a}_{n,s+q}}}},\]
\[k,n=1,2,3,\ldots ,\quad p,q=0,1,\ldots .\]
The least favorable stochastic sequence in the class $\mathbf{Y}$ for the optimal estimate of the functional $A\zeta $ is a one-sided moving average sequence of the form
\[{{\vec{\zeta }}_{j}}=\sum\nolimits_{s=-\infty }^{j}{g(j-s)\vec{\varepsilon }(s)},\]
where $g=(g(p))_{p=0}^{\infty }$ is the eigenvector, that corresponds to $\nu^2$, constructed from matrices $g(s)=\{{{g}_{km}}(s)\}_{k=\overline{1,\infty }}^{m=\overline{1,M}}$  and  determined by the condition $\|g{{\|}^{2}}=\sum\nolimits_{p=0}^{\infty }{\|g(p){{\|}^{2}}}=P$, and $\vec{\varepsilon }(s)=$ $\{{{\varepsilon }_{m}}(s)\}_{m=1}^{M}$ is a vector stationary stochastic sequence with orthogonal values.
\end{thm}

   Proof:  Upper bound. Let us consider the norm approximation of the self-adjoint continuous [17] operator $Q=A{{A}^{*}},$ $A=\{A(p,q)\}_{p,q=0}^{\infty }=\{{{\vec{a}}_{p+q}}\}_{p,q=0}^{\infty }$ in ${{\ell }_{2}}$ by  the sequence of operators  ${{Q}_{N}}={{A}_{N}}A_{N}^{*}$ determined in accordance with Theorem \ref{thm3.1}.
   Then we have
\[\underset{\hat{A}\in \Lambda }{\mathop{\min }}\,\underset{\zeta \in \mathbf{Y}}{\mathop{\max }}\,\Delta (\zeta ,\hat{A})=\underset{N\to \infty }{\mathop{\lim }}\,(\underset{{{{\hat{A}}}_{N}}\in \Lambda }{\mathop{\min }}\,\underset{\zeta \in \mathbf{Y}}{\mathop{\max }}\,\Delta (\zeta ,{{\hat{A}}_{N}}))=\]
\begin{equation}\label{eq13}
=P\underset{N\to \infty }{\mathop{\lim }}\,\nu _{N}^{2}=P\cdot {{\nu }^{2}}.
\end{equation}

   Lower bound. Solution of the extremal problem $$|(Q\tilde{g},\tilde{g})|\to \max,\|\tilde{g}\|=1$$ is an eigenvector,
   that corresponds to the greatest eigenvalue [17]
  $${{\nu }^{2}}=\underset{\|\tilde{g}\|\le 1}{\mathop{\max }}\,|(Q\tilde{g},\tilde{g})|=\|Q\|.$$
  Thus, the following inequalities hold true
\[\underset{\zeta \in {{\mathbf{Y}}_{\mathbf{R}}}}{\mathop{\max }}\,\underset{\hat{A}\in \Lambda }{\mathop{\min }}\,\Delta (\zeta ,\hat{A})=P\cdot \underset{\|\tilde{g}\|\le 1}{\mathop{\max }}\,|(Q\tilde{g},\tilde{g})|=\]
\begin{equation}\label{eq14}
=P\cdot {{\nu }^{2}}\le \underset{\zeta \in \mathbf{Y}}{\mathop{\max }}\,\underset{\hat{A}\in \Lambda }{\mathop{\min }}\,\Delta (\zeta ,\hat{A}).
\end{equation}
From relations \eqref{eq13} and \eqref{eq14}  equality \eqref{eq12} follows. Theorem is proved.

\section{\Large{Conclusions}}

   Formulas for calculating the maximum values of the mean square errors of optimal linear estimation of the functionals
   $A\zeta $ and ${{A}_{N}}\zeta $ that depend on the unknown values of periodically correlated stochastic process $\zeta (t)$
   from the class $\mathbf{Y}$ are proposed. The estimation is based on observations of the process $\zeta (t)$ for $t<0$.
   It is shown that the maximum value of the error of estimate of the functional $A\zeta$ (as well as ${{A}_{N}}\zeta $) give the one-sided moving average stationary sequence which corresponds to the PC process $\zeta (t)$.

\end{document}